\documentclass[10pt]{article}
\usepackage[fleqn]{amsmath}
\usepackage{amsmath,amsthm,amsfonts,amssymb,comment}

\allowdisplaybreaks

\newcommand{\p}{\partial}
\newcommand\pd[2]{\frac{\partial#1}{\partial#2}}

\numberwithin{equation}{section}

\begin{document}

\title{\Large\bf Moufang symmetry I.\\Generalized Lie and Maurer-Cartan equations}
\date{}

\author{Eugen Paal}
\maketitle
\thispagestyle{empty}

\begin{abstract}
The continuous Moufang loops are characterized as the algebraic systems where the associativity law is perturbed minimally. The minimal perturbation of associativity is characterized by the first-order partial differential equations, which in a natural way generalize the Lie and Maurer-Cartan equations from the theory of Lie groups. 
\end{abstract}


\section{Introduction and outline of the paper}

In 1972, Sabinin stated \cite{Sab72a,Sab72b} equivalence of algebraic loops and homogeneous spaces. Following his construction, one may characterize loops as  \textit{groups with hidden associativity}. Due to perturbation  of   associativity, the particular loops are too difficult to describe. In such a situation, the Moufang loops are distinguished by their \textit{minimal} perturbation of associativity. Thus, one may expect that a Moufang loop represents a minimally perturbed symmetry that is concisely referred as the \textit{Moufang symmetry}.

In this paper, the minimal perturbation of associativity in the continuous Mofang loops is characterized by the first-order partial differential equations, which in a natural way generalize the Lie and Maurer-Cartan equations from the theory of Lie groups. 

\section{Moufang loops}

A \textit{Moufang loop} \cite{RM} (see also \cite{Bruck,Bel,HP}) is a set $G$ with a binary operation  (multiplication) 
$\cdot: G\times G\to G$,
denoted also by juxtaposition, so that the following three axioms are satisfied.
\begin{enumerate}
\itemsep-2pt
\item[1)] 
In the equation $gh=k$, the knowledge of any two of $g,h,k\in G$ specifies the third one \textit{uniquely}.
\item[2)] 
There is a distinguished element $e\in G$ with the property $eg=ge=g$ for all $g\in G$.
\item[3)] 
The \textit{Moufang identities} hold in $G$, 
\begin{align}
\label{m}
   g(h\cdot gk) = (gh\cdot g)k,\quad
   (kg\cdot h)g = k(g\cdot hg),\quad
   (gh)(kg) = g(hk\cdot g)
\end{align}
\end{enumerate}
It must be noted that the Moufang identities (\ref{m}) are equivalent (see e.g \cite{Bruck,Bel}), which means that each of the above identities may be considered as a defining identity of the Moufang loop.

Recall that a set with a binary operation is called a \textit{groupoid}. A groupoid $G$ with axiom 1) is called a \textit{quasigroup}. If axioms 1) and 2) are satisfied, the groupoid (quasigroup) $G$ is called a \textit{loop}. The element $e$ in axiom 2) is called the \textit{unit} (element) of the (Moufang) loop $G$.

In a (Moufang) loop, the multiplication need not be neither associative nor commutative. The \textit{associativity} and \textit{commutativity} laws read, respectively,
\begin{equation*}
g\cdot hk=gh\cdot k,\quad gh=hg, \quad  \forall\,  g,h,k\in G
\end{equation*}
The most familiar kind of loops are the ones with the \textit{associative} law, and these are called \textit{groups}.   A  loop $G$ is called \textit{commutative} if the commutativity law holds in $G$, and (only) the commutative associative loops are said to be \textit{Abelian}.

A characteristic property of the Moufang loops is their \textit{diassociativity}: in a Moufang loop  every two elements generate an associative subloop (group) \cite{RM}. In particular, from this it follows that
\begin{align*} 
g\cdot gh=g^{2}h,\quad 
hg\cdot g=hg^{2},\quad 
gh\cdot g=g\cdot hg,\quad \forall \, g,h\in G
\end{align*}
Here, the first two identities are called the left and right \textit{alternativity}, respectively, and the third one is said to be \textit{flexibility}. Note that these identities follow from the Moufang identities. Due to flexibility, the Moufang identity can be rewritten in a nice symmetric  form
\begin{equation*} 
gh\cdot kg=g\cdot hk\cdot g,\quad \forall\, g,h,k\in G 
\end{equation*}

The unique solution of the equation $xg=e$ ($gx=e$) is called the left (right) \textit{inverse} element of $g\in G$ and is denoted as $g^{-1}_{R}$ ($g^{-1}_{L}$). It also follows from the diassociativity of the Moufang loop that 
\begin{align*}
g^{-1}_{R}=g^{-1}_{L}\doteq  g^{-1},\quad
g^{-1}\cdot gh = hg\cdot g^{-1}=h,\quad 
\left(g^{-1}\right)^{-1}=g,\quad 
(gh)^{-1}=h^{-1}g^{-1}
\end{align*}
for all $g,h \in G$.

The \textit{smallest} non-associative Moufang loop was described by Chein and Pflugfelder in \cite{CP71}, its order is 12. One can find an elaborated list of the low-order non-associative Moufang loops in \cite{Chein78}.

\section{Analytic Moufang loops and tangent Malcev algebras}

A Moufang loop $G$ is said to be \textit{analytic} \cite{Mal} if $G$ is a finite dimensional real, analytic manifold so that both the Moufang loop operation $G\times G\to G$: ($g,h$) $\mapsto gh$ and the inversion map $G\to G$: $g\mapsto g^{-1}$ are analytic ones. The dimension of $G$ will be denoted as $\dim G\doteq  r$. 

In this paper, the global aspects are not considered and all considerations are local, i.e in a fixed chart $U_e$ of the unit $e\in U_e\subset G$ where all subsequent constructions are expected to be well defined. The triple $(G,e,U_e)$ is called a \text{local} analytic Moufang loop. By keeping this in mind we often call $G$  the \textit{local} analytic Moufang loop as well.

The local coordinates of $g\in G$ are denoted (in selected fixed chart $U_e\subset G$ of the unit $e\in G$) by $g^{1},\dots,g^{r}$ so that the local coordinates of the unit $e$ are supposed to be zero: $e^{i}=0$, $i=1,\dots,r$. One has the evident initial conditions
\begin{equation*}
(ge)^{i}=(eg)^{i}=g^{i},\quad i=1,\dots,r
\end{equation*}
As in the case of the Lie groups \cite{Pontr}, we can use the Taylor expansions
%
\begin{align*}
(gh)^{i}
&=h^{i}+u^{i}_{j}(h)g^{j}+\frac{1}{2!}u^{i}_{jk}(h)g^{j}h^{k}\cdots\\
&=g^{i}+v^{i}_{j}(g)h^{j}+\frac{1}{2!}v^{i}_{jk}(g)g^{j}h^{k}\cdots\\
&=g^{i}+h^{i}+a^{i}_{jk}g^{j}h^{k}
+b^{i}_{jkl}g^{j}g^{k}h^{l}
+d^{i}_{jkl}g^{j}h^{k}h^{l}+\cdots
\end{align*}
%
to introduce the \textit{auxiliary functions} $u^{i}_{j}, v^{i}_{j}, u^{i}_{jk}, v^{i}_{jk}$ and the \textit{structure constants}
\begin{equation*} 
c^{i}_{jk}\doteq  a^{i}_{jk}-a^{i}_{kj}=-c^{i}_{kj}   
\end{equation*}
Due to such a parametrization, one has the initial conditions
\begin{align*}
u^{i}_{j}(e)=\delta^{i}_{j}=v^{i}_{j}(e),\quad
u^{i}_{jk}(e)=0=v^{i}_{jk}(e),\quad
\end{align*}
and symmetry with respect to the lower indices 
\begin{align*}
u^{i}_{jk}=u^{i}_{kj},\quad
v^{i}_{jk}=v^{i}_{kj}
\end{align*}
Also, it follows from axiom 1) of the Moufang loop that 
\begin{equation*}
\det[u^{i}_{j}]\neq0\neq\det[v^{i}_{j}]
\end{equation*}
The \textit{tangent algebra} of $G$ can be defined similarly to the tangent algebra of the Lie group 
\cite{Pontr}. Geometrically, this algebra is the tangent space $T_{e}(G)$ of $G$ at $e$. 
The product of $x,y\in T_{e}(G)$ will be denoted by $[x,y]\in T_{e}(G)$. In coordinate form,
\begin{equation*} 
[x,y]^{i}\doteq  c^{i}_{jk}x^{j}y^{k}
            =-[y,x]^{i},\quad i=1,\dots,r 
\end{equation*}
The tangent algebra of $G$ will be denoted by  $\Gamma_G\doteq (T_{e}(G),[\cdot,\cdot])$. The latter need not be a Lie algebra. In other words, there may exist a triple $x,y,z\in T_{e}(G)$, such that the Jacobi identity fails,
\begin{equation*} 
J(x,y,z)\doteq  [x,[y,z]]+[y,[z,x]]+[z,[x,y]]\ne0 
\end{equation*}
Instead, for all $x,y,z\in T_{e}(G)$, one has \cite{Mal} a more general identity
\begin{equation*} 
[[x,y],[z,x]]+[[[x,y],z],x]+[[[y,z],x],x]+[[[z,x],x],y]=0
\end{equation*}
called the \textit{Malcev identity}. The tangent algebra $\Gamma_G$  is hence said to be the \textit{Malcev algebra}. The Malcev identity concisely reads \cite{Sagle}
\begin{equation*}
[J(x,y,x),x]=J(x,y,[x,z])
\end{equation*}
from which it can be easily seen that every Lie algebra is a Malcev algebra as well. It has been shown in  \cite{Kuzmin} that every  finite-dimensional real Malcev algebra is the tangent algebra of some local analytic Moufang loop.

The \textit{smallest} non-associative real analytic Moufang loop and its Malcev algebra were described by Akivis in \cite{Akivis77}, their dimension is 4.

\section{Associators}

In a (Moufang) loop, due to non-associativity, the elements $g\cdot hk$ and $gh\cdot k$ need not coincide. 
The non-associativity of $G$ can be measured by the  formal non-invariant functions
\begin{equation*} 
A^{i}(g,h,k)\doteq  (g\cdot hk)^{i}-(gh\cdot k)^{i},\quad i=1,\dots,r
\end{equation*}
which are called \textit{associators} of $G$. One has the evident initial conditions
\begin{equation*} 
A^{i}(e,h,k)=A^{i}(g,e,k)=A^{i}(g,h,e)=0,\quad i=1,\dots,r
\end{equation*}
The associators $A^{i}$ are considered as generating functions in the following sense. At first, define the 
\textit{first-order} associators $l^{i}_{j}$, $r^{i}_{j}$, $m^{i}_{j}$ by
\begin{align*}
A^{i}(g,h,k)
&\doteq  l^{i}_{j}(h,k)g^{j}+O(g^{2})\\
&\doteq  r^{i}_{j}(g,h)k^{j}+O(k^{2})\\
&\doteq  m^{i}_{j}(g,k)h^{j}+O(h^{2})
\end{align*} 
These can be easily calculated and the result reads
%
\begin{align*}
l^{i}_{j}(g,h)  &=-u^{s}_{j}(g)\pd{(gh)^{i}}{g^{s}}+u^{i}_{j}(gh)\\
r^{i}_{j}(g,h)  &=-v^{i}_{j}(gh)+v^{s}_{j}(h)\pd{(gh)^{i}}{h^{s}}\\
m^{i}_{j}(g,h)&=-v^{s}_{j}(g)\pd{(gh)^{i}}{g^{s}}+u^{s}_{j}(h)\pd{(gh)^{i}}{h^{s}}
\end{align*}
%
Next one can check the initial conditions
\begin{align*}
l^{i}_{j}(e,g)=r^{i}_{j}(e,g)=m^{i}_{j}(e,g)
=0=
l^{i}_{j}(g,e)=r^{i}_{j}(g,e)=m^{i}_{j}(g,e)
\end{align*}
and define the \textit{second-order} associators
$l^{i}_{jk}$, $\tilde l^{i}_{jk}$, $m^{i}_{jk}$, $\tilde m^{i}_{jk}$, $r^{i}_{jk}$, $\tilde r^{i}_{jk}$ by
\begin{align*} 
l^{i}_{j}(g,h)&\doteq  l^{i}_{jk}(h)g^{k}+O(g^{2})\\
              &\doteq  \tilde l^{i}_{jk}(g)h^{k}+O(h^{2})\\
r^{i}_{j}(g,h)&\doteq  r^{i}_{jk}(h)g^{k}+O(g^{2})\\
             &\doteq  \tilde r^{i}_{jk}(g)h^{k}+O(h^{2})\\
m^{i}_{j}(g,h)&\doteq  \tilde m^{i}_{jk}(h)g^{k}+O(g^{2})\\
              &\doteq  m^{i}_{jk}(g)h^{k}+O(h^{2})
\end{align*}
By calculating, the result reads
%
\begin{align*} 
&l^{i}_{jk}(g)=\tilde m^{i}_{kj}(g)                           
             =-u^{i}_{kj}(g)-a^{s}_{jk}u^{i}_{s}(g)
               +u^{s}_{k}(g)\pd{u^{i}_{j}(g)}{g^{s}}\\
&r^{i}_{jk}(g)=\tilde l^{i}_{kj}(g)                              
             = v^{s}_{j}(g)\pd{u^{i}_{k}(g)}{g^{s}}
               -u^{s}_{k}(g)\pd{v^{i}_{j}(g)}{g^{s}}  \\
&m^{i}_{jk}(g)
=\tilde r^{i}_{kj}(g)                          
             =v^{i}_{jk}(g)+a^{s}_{jk}v^{i}_{s}(g)
               -v^{s}_{j}(g)\pd{v^{i}_{k}(g)}{g^{s}}                
\end{align*}
Finally one has the initial conditions
\begin{align*}
l^{i}_{jk}(e)=r^{i}_{jk}(e)=m^{i}_{jk}(e)
=0=
\tilde l^{i}_{jk}(e)=\tilde r^{i}_{jk}(e)=\tilde m^{i}_{jk}(e)
\end{align*}
and define the \textit{third-order} associators
         $l^{i}_{jkl}$, $\tilde l^{i}_{jkl}$, 
         $m^{i}_{jkl}$, $\tilde m^{i}_{jkl}$,
         $r^{i}_{jkl}$, $\tilde r^{i}_{jkl}$ by
\begin{align*}
&l^{i}_{jk}(g)\doteq  l^{i}_{jkl}g^{l}+O(g^{2}),\quad\hspace{13pt}
\tilde l^{i}_{jk}(g)\doteq  \tilde l^{i}_{jkl}g^{l}+O(g^{2})\\
&r^{i}_{jk}(g)\doteq  r^{i}_{jkl}g^{l}+O(g^{2}),\quad \hspace{9pt}
\tilde r^{i}_{jk}(g)\doteq  \tilde r^{i}_{jkl}g^{l}+O(g^{2})\\
&m^{i}_{jk}(g)\doteq  m^{i}_{jkl}g^{l}+O(g^{2}),\quad
\tilde m^{i}_{jk}(g)\doteq   \tilde m^{i}_{jkl}g^{l}+O(g^{2})
\end{align*}
By calculating, one gets
\begin{subequations}
\begin{align*}
l^{i}_{jkl}= m^{i}_{klj}=r^{i}_{ljk}
            = \tilde l^{i}_{jlk}=\tilde r^{i}_{lkj}=\tilde m^{i}_{kjl}
           = a^{i}_{js}a^{s}_{kl}-a^{s}_{jk}a^{i}_{sl}
             +2\left(d^{i}_{jkl}-b^{i}_{jkl}\right)
\end{align*}
\end{subequations}
Thus, all third-order associators can be obtained from $l^{i}_{jkl}$ via permutations of the lower indices. One can check the \textit{Akivis identity} \cite{Akivis}
\begin{subequations}
\label{Akivis}
\begin{align*}
l^{i}_{jkl}+r^{i}_{jkl}+m^{i}_{jkl}-\tilde l^{i}_{jkl}-\tilde r^{i}_{jkl}-\tilde m^{i}_{jkl}
&=l^{i}_{jkl}+l^{i}_{klj}+l^{i}_{ljk}-l^{i}_{jlk}-l^{i}_{lkj}-l^{i}_{kjl}\\ 
&=c^{i}_{js}c^{s}_{kl}+c^{i}_{ks}c^{s}_{lj}+c^{i}_{ls}c^{s}_{jk}
\end{align*}
\end{subequations}
For $x,y,z\in T_e(G)$ define \cite{Akivis} their trilinear product $(x,y,z)\in T_e(G)$ by
\begin{equation*} 
(x,y,z)^{i}\doteq  l^{i}_{jkl}x^{j}y^{k}z^{l},\quad i=1\dots,r
\end{equation*}
Then the Akivis identity reads \cite{Akivis}
\begin{align*}
J(x,y,z)=(x,y,z)+(y,z,x)+(z,x,y)-(x,z,y)-(z,y,x)-(y,x,z)
\end{align*}

\section{Generalized Lie equations}

Differentiate the Moufang identities (\ref{m}) with respect to $g^{j}$ at $g=e$. Then, suitable re-denoting variables we obtain for the first-order associators the constraints
%
\label{m_d1}
\begin{align*} 
l^{i}_{j}(g,h)=-m^{i}_{j}(g,h),\quad
r^{i}_{j}(g,h)=-m^{i}_{j}(g,h),\quad
l^{i}_{j}(g,h)=r^{i}_{j}(g,h)
\end{align*}
%
Thus, the Moufang identities give rise to the first-order constraints
\begin{equation}
\label{min1}
l^{i}_{j}(g,h)=r^{i}_{j}(g,h)=-m^{i}_{j}(g,h)
\end{equation}
The latter in a natural way generalize the  Lie equations. If $G$ is associative, then one has  the Lie equations (see e.g \cite{Pontr}):
\begin{align*}
l^{i}_{j}(g,h)=r^{i}_{j}(g,h)=-m^{i}_{j}(g,h)=0
\end{align*}
By comparing the latter with (\ref{min1}) one can say that in the variety of continuous non-associative loops the local analytic Moufang loops have the property that their associativity is perturbed  \textit{minimally}. Constraints (\ref{min1}) are hence called the 
\textit{first-order minimality conditions}.

Define the new auxiliary functions $w^{i}_{j}$ by
\begin{equation}
\label{def_w}
u^{i}_{j}(g)+v^{i}_{j}(g)+w^{i}_{j}(g)=0
\end{equation}
Then constraints (\ref{min1}) can be rewritten as the differential equations of the 1st order, called the 
\textit{generalized Lie equations}
\begin{align*}
w^{s}_{j}(g)\pd{(gh)^{i}}{g^{s}}+u^{s}_{j}(h)\pd{(gh)^{i}}{h^{s}}+u^{i}_{j}(gh)&=0\\
v^{s}_{j}(g)\pd{(gh)^{i}}{g^{s}}+w^{s}_{j}(h)\pd{(gh)^{i}}{h^{s}}+v^{i}_{j}(gh)&=0\\
u^{s}_{j}(g)\pd{(gh)^{i}}{g^{s}}+v^{s}_{j}(h)\pd{(gh)^{i}}{h^{s}}+w^{i}_{j}(gh)&=0
\end{align*}
In view of (\ref{def_w}) these equations are linearly dependent.

\section{Generalized Maurer-Cartan equations}

Differentiate constraints (\ref{min1})  with respect to $g^{k}$ and $h^{k}$ at $g=e=h$, respectively. Then, re-denoting variables one obtains the second-order constraints
\begin{equation}
\label{m_d2}
l^{i}_{jk}(g)=r^{i}_{jk}(g)=m^{i}_{jk}(g)=-m^{i}_{kj}(g) 
\end{equation}
Again, if $G$ is associative, i.e a local Lie group, we have
\begin{equation*}
l^{i}_{jk}(g)=r^{i}_{jk}(g)=m^{i}_{jk}(g)=-m^{i}_{kj}(g)=0 
\end{equation*}
Constraints (\ref{m_d2}) are called the  \textit{second-order minimality conditions} for $G$. Lets consider the latter more closely.

It follows from the skew-symmetry $l^{i}_{jk}=-l^{i}_{kj}$ and $m^{i}_{jk}=-m^{i}_{kj}$, respectively, that
%
\begin{align*}
2u^{i}_{jk}
=u^{s}_{k}\pd{u^{i}_{j}}{g^{s}}+u^{s}_{j}\pd{u^{i}_{k}}{g^{s}}
            -\left(a^{i}_{jk}+a^{i}_{kj}\right)u^{s}_{j},\quad
2v^{i}_{jk}
=v^{s}_{k}\pd{v^{i}_{j}}{g^{s}}+v^{s}_{j}\pd{v^{i}_{k}}{g^{s}}
            -\left(a^{i}_{jk}+a^{i}_{kj}\right)v^{s}_{j}
\end{align*}
%
By using $u^{i}_{jk}$ and $v^{i}_{jk}$  one can see that 
\begin{align*}
u^{s}_{k}\pd{u^{i}_{j}}{g^{s}}-u^{s}_{j}\pd{u^{i}_{k}}{g^{s}}
             =c^{s}_{jk}u^{i}_{s}+2l^{i}_{jk},\quad
v^{s}_{k}\pd{v^{i}_{j}}{g^{s}}-v^{s}_{j}\pd{v^{i}_{k}}{g^{s}}
             =c^{s}_{kj}v^{i}_{s}+2m^{i}_{jk}
\end{align*}
Now use $l^{i}_{jk}=r^{i}_{jk}$, $m^{i}_{jk}=-r^{i}_{kj}$ and $r^{i}_{jk}$ to obtain
the differential equations for the first order auxiliary functions $u^{i}_{j}$ and  $v^{i}_{j}$, 
\begin{subequations}
\label{gen-m-c}
\begin{align}
u^{s}_{k}\pd{u^{i}_{j}}{g^{s}}-u^{s}_{j}\pd{u^{i}_{k}}{g^{s}}
&=c^{s}_{jk}u^{i}_{s} +2\left(v^{s}_{j}\pd{u^{i}_{k}}{g^{s}}
                           -u^{s}_{k}\pd{v^{i}_{j}}{g^{s}}\right)\\
v^{s}_{k}\pd{v^{i}_{j}}{g^{s}}-v^{s}_{j}\pd{v^{i}_{k}}{g^{s}}
&=c^{s}_{kj}v^{i}_{s}+2\left(u^{s}_{j}\pd{v^{i}_{k}}{g^{s}}
                           -v^{s}_{k}\pd{u^{i}_{j}}{g^{s}}\right)
\end{align}
\end{subequations}
called the \textit{generalized Maurer-Cartan equations} of $G$, thees  generalize the Maurer-Cartan equations (see e.g \cite{Pontr}) in a minimal way.

Finally, note that the constraints $r^{i}_{jk}=-r^{i}_{kj}$ read
\begin{equation}
\label{l-r}
v^{s}_{j}\pd{u^{i}_{k}}{g^{s}}-u^{s}_{k}\pd{v^{i}_{j}}{g^{s}}
  =u^{s}_{j}\pd{v^{i}_{k}}{g^{s}}-v^{s}_{k}\pd{u^{i}_{j}}{g^{s}}
\end{equation}

The generalized Maurer-Cartan differential equations can be rewritten more concisely. 
For $x\in T_e(G)$ introduce the \textit{infinitesimal translations} as the local vector fields in the vicinity of the unit $e\in G$,
\begin{equation*}
L_x \doteq L_x(g) \doteq x^j u^{i}_j(g) \frac{\p}{\p g^{i}}, \quad 
R_x \doteq R_x(g) \doteq x^j v^{i}_j(g) \frac{\p}{\p g^{i}} \quad \in T_g(G)
\end{equation*}
Define the Lie bracketing $[A,B]$ of the (local) vector fields $A$ and $B$ (on $G$) in the usual way $[A,B]\doteq  AB-BA$. Then,  the generalized Maurer-Cartan  equations (\ref{gen-m-c}a,b) and differential equations (\ref{l-r}) can be rewritten, respectively, as the commutation relations (CR) called a \textit{birepresentation} \cite{Paal93} of the tangent Malcev algebra $\Gamma_G$ as follows:
%
\begin{align*}
\left[L_{x},L_{y}\right]=L_{[x,y]}-2\left[L_{x},R_{y}\right] ,\quad 
\left[R_{x},R_{y}\right]=R_{[y,x]}-2\left[R_{x},L_{y}\right],\quad 
\left[L_{x},R_{y}\right]=\left[R_{x},L_{y}\right]
\end{align*}
%
Here, one can easily see relation to the structure equations of \textit{alternative} algebras \cite{Schafer66}. 
In the case of associativity of $G$, the (generalized) Maurer-Cartan equations read
\begin{equation*}
2\left[L_{x},R_{y}\right]
=L_{[x,y]}-\left[L_{x},L_{y}\right]  
=R_{[y,x]}-\left[R_{x},R_{y}\right]
=2\left[R_{x},L_{y}\right]=0 
\end{equation*}

\section{Akivis identity}

Differentiate (\ref{m_d2}) with respect to $g^l$ ($l=1,\ldots,n$) at $g=e$. Then one gets the 3rd order) constraints, called the \textit{third-order minimality conditions}, 
\begin{equation*}
l^{i}_{jkl}=l^{i}_{klj}=-l^{i}_{lkj}=r^{i}_{jkl}=m^{i}_{jkl}
=-\tilde{l}^{i}_{jkl}=-\tilde{r}^{i}_{jkl}=-\tilde{m}^{i}_{jkl}
\end{equation*}
Thus the third-order associators of a local analytic Moufang loop have the property of their total anti-symmetry with respect of the lower indexes. From this it follows that for a continuous  Moufang loop $G$ the Akivis identity reads 
\begin{align*}
J(x,y,z)=6(x,y,z),\quad \forall \, x,y,z \in T_e(G)
\end{align*}

\section*{Acknowledgements}

The research was in part supported by the Estonian Research Council (Grant ETF9038) and Astralgo Tallinn. Author also partially used, modified  and extended material from \cite{Paal87}.

\par\noindent
\footnotesize{Tallinn University of Technology, Ehitajate tee 5, 19086 Tallinn, Estonia.\\
Astralgo, Tallinn, Estonia}
\end{document}